\documentclass{article}
\usepackage{amssymb}
\usepackage{amsmath}
\usepackage{graphicx, anysize}
\usepackage{multirow}
\usepackage{epstopdf}
\usepackage{natbib}
\usepackage{color}
\begin{document}

\newtheorem{axiom}{Axiom}[section]
\newtheorem{definition}{Definition}[section]
\newtheorem{statement}{Statement}[section]
\newtheorem{comment}{Comment}[section]
\newtheorem{convention}{Convention}[section]
\newtheorem{proposition}{Proposition}[section]
\newtheorem{lemma}{Lemma}[section]
\newtheorem{theorem}{Theorem}[section]
\newtheorem{corollary}{Corollary}[section]
\newtheorem{remark}{Remark}[section]

\newenvironment{proof}%
{\noindent}{\hfill$\square$}%

\setlength{\baselineskip}{12pt}

\begin{center} \large \sc On Asymptotic Normality of the Local Polynomial Regression Estimator with Stochastic Bandwidths\footnote{We thank an anonymous referee for very useful suggestions.  We also thank Juan Carlos Escanciano, Yanqin Fan and Jeff Racine for helpful comments.  We are particularly grateful to Yanqin Fan for bringing to our attention the work of Ziegler (2004). } \normalsize \rm  \\[.5in]

\end{center}
\begin{center}
\begin{tabular}{lcl}
\multicolumn{3}{c}{\large \sc Carlos Martins-Filho}\\[.2in] 
Department of Economics &  &IFPRI \\  
University of Colorado &  & 2033 K Street NW  \\
Boulder, CO 80309-0256, USA& \&&Washington, DC 20006-1002, USA\\ 
email: carlos.martins@colorado.edu& & email: c.martins-filho@cgiar.org\\
Voice: + 1 303 492 4599 & &Voice: + 1 202 862 8144\\[.5pt]
\end{tabular}\\[.2in]

and\\[.2in]

\begin{tabular}{l}
\multicolumn{1}{c}{\large \sc Paulo Saraiva}\\[.2in]
Department of Economics \\  
University of Colorado \\
Boulder, CO 80309-0256, USA\\ 
email: paulo.saraiva@colorado.edu\\
Voice: +1 541 740 5209 \\[.4in]
\end{tabular}

October, 2010\\[.4in]
\end{center}

\noindent\bf Abstract. \rm Nonparametric density and regression estimators commonly depend on a bandwidth. The asymptotic properties of these estimators have been widely studied when bandwidths are nonstochastic. In practice, however, in order to improve finite sample performance of these estimators, bandwidths are selected by data driven methods, such as cross-validation or plug-in procedures. As a result nonparametric estimators are usually constructed using stochastic bandwidths. In this paper we establish the asymptotic equivalence in probability of local polynomial regression estimators under stochastic and nonstochastic bandwidths. Our result extends previous work by \cite{Boente1995} and \cite{Ziegler2004}.\\[.1in]



\clearpage

\definecolor{red}{rgb}{1,0,0}
\setlength{\baselineskip}{24pt}
\pagestyle{plain}
\setcounter{page}{1}
\setcounter{footnote}{0}
\section{Introduction}\label{intro}
Currently there exist several papers that establish the asymptotic properties of kernel based nonparametric estimators.  For the case of density estimation, \cite{Parzen1962}, \cite{Robinson1983} and \cite{Bosq1998} establish the asymptotic normality of Rosenblatt's density estimator under independent and identically distributed (IID)  and stationary strong mixing data generating processes.  For the case of regression, \cite{Fan1992}, \cite{Masry1997} and \cite{Martinsfilho2009} establish asymptotic normality of local polynomial estimators under IID, stationary and nonstationary strong mixing processes.  All of these asymptotic approximations are obtained for a sequence of nonstochastic bandwidths $0<h_n \rightarrow 0$ as the sample size $n\rightarrow \infty$.

In practice, to improve estimators' finite sample performance, bandwidths are normally selected using data-driven methods (\cite{Ruppert1995,  Xia2002}). As such, bandwidths are in practical use generally stochastic. Therefore, it is desirable to obtain the aforementioned asymptotic results when $h_n$ is data dependent and consequently stochastic.

There have been previous efforts in establishing asymptotic properties of nonparametric estimators constructed with stochastic bandwidths. Consider, for example, the local polynomial regression estimator proposed by \cite{Fan1992}.  \cite{Dony2006} prove that such estimator, when constructed with a stochastic bandwidth, is uniformly consistent. More precisely, suppose $\{(Y_t,X_t)\}_{t=1}^n$ is a sequence of random vectors in $\mathbb{R}^2$  with regression function $m(x)=E(Y_t|X_t=x)$ for all $t$.  The local polynomial regression estimator of order $p$ is defined by $m_{LP}(x;h_n)\equiv\hat{b}_{n0}(x;h_n)$ where
$$(\hat{b}_{n0}(x;h_n),\ldots,\hat{b}_{np}(x;h_n))=\underset{b_0,\ldots,b_p}{\mbox{argmin}}\sum_{t=1}^n\left(Y_t-\sum_{j=0}^pb_j(X_t-x)^j \right)^2K\left(\frac{X_t-x}{h_n} \right)$$
and $K:\mathbb{R}\rightarrow\mathbb{R}$ is a kernel function. If the sequence $\{(Y_t,X_t)\}_{t=1}^n$ is IID, then it follows from \cite{Dony2006} that 
$$\limsup_{n\rightarrow\infty}\sup_{h_n\in[a_n, b_n]}\frac{\sqrt{nh_n}\sup_{x\in G}|m_{LP}(x;h_n)-m(x)|}{\sqrt{|\log h_n|\lor \log\log n}}=O_{a.s.}(1)$$
where $(a_n,b_n)$ is a nonstochastic sequence such that $0\leq a_n<b_n\rightarrow 0$ as $n\rightarrow\infty$, $G$ is a compact set in $\mathbb{R}$ and $|\log h_n|\lor \log\log n=\max\{|\log h_n|, \log\log n\}$. 
If there exists a stochastic bandwidth $\hat{h}_n$ such that $\frac{\hat{h}_n}{h_n}-1=o_p(1)$ and we define $a_n=rh_n$ and $b_n=sh_n$ with $0<r<1<s$. Then it follows that
$$\sup_{x\in G}|m_{LP}(x;\hat{h}_n)-m(x)|=o_p(1).$$

When $p=1$ and the sequence $\{(Y_t,X_t)\}_{t=1}^n$ is IID, if $\hat{h}_n$ is obtained by a cross validation procedure, \cite{Li2004} show that
$$\sqrt{n\hat{h}_n}\left(m_{LP}(x;\hat{h}_n)-m(x)-\frac{\hat{h}^2_n}{2}m^{(2)}(x)\int K(u)u^2du\right)\stackrel{d}{\rightarrow}N\left(0,\frac{\sigma^2(x)}{f_X(x)}\int K^2(u)du\right)$$
where $X$ is a random variable that has the same distribution of $X_t$, $f_X$ is the density function of $X$ and $\sigma^2(x)=Var(Y_t|X_t=x)$. \cite{Xia2002} establish that, if $\hat{h}_n$ is obtained through cross validation, $\frac{\hat{h}_n}{h_n}-1=o_p(1)$ for strong mixing and strictly stationary sequences $\{(Y_t,X_t)\}_{t=1}^n$.

When $p=0$, the case of a Nadaraya-Watson regression estimator $m_{NW}(x;h_n)$, and the sequence $\{(Y_t,X_t)\}_{t=1}^n$ is a strictly stationary strong mixing random process, \cite{Boente1995} show that if $\frac{\hat{h}_n}{h_n}-1=o_p(1)$, then
$$\sqrt{n\hat{h}_n}\left(m_{NW}(x;\hat{h}_n)-E\left(m_{NW}(x;\hat{h}_n)\left|\vec{X}\right.\right)\right)\stackrel{d}{\rightarrow}N\left(0,\frac{\sigma^2(x)}{f_X(x)}\int K^2(u)du\right)$$
where $\vec{X}'=(X_1,\ldots,X_n)$.  Since independent processes are strong mixing, their result encompasses the case where $\{(Y_t,X_t)\}_{t=1}^n$ is IID, which is treated in the otherwise broader paper by \cite{Ziegler2004}.

In this paper we expand the result of \cite{Boente1995} by obtaining that local polynomial estimators for the regression and derivatives of orders $j=1,\ldots,p$ constructed with a stochastic bandwidth $\hat{h}_n$ are asymptotically normal. We do this for processes that are strong mixing and strictly stationary. Our proofs build and expand on those of \cite{Boente1995} and \cite{Masry1997}.

\section{Preliminary Results and Assumptions}

Define the vector $b_n(x;h)=(\hat{b}_{n0}(x;h),\ldots,\hat{b}_{np}(x;h))'$ and the diagonal matrix $H_n=diag\{h_n^j\}_{j=0}^p$.  Given that \cite{Masry1997} have established the asymptotic normality of $\sqrt{nh_n}\left(H_nb_n(x;h_n)-E(H_nb_n(x;h_n)|\vec{X})\right)$, it suffices for our purpose to show that
\begin{equation}\label{eq:obj}
\sqrt{nh_n}\left(H_n(h_n)b_n(x;h_n)-E(H_nb_n(x;h_n)|\vec{X})\right)-\sqrt{n\hat{h}_n}\left(\hat{H}_nb_n(x;\hat{h}_n)-E(\hat{H}_nb_n(x;\hat{h}_n)|\vec{X})\right)=o_p(1),
\end{equation}
where $\hat{h}_n$ is a bandwidth that satisfies $\frac{\hat{h}_n}{h_n}-1=o_p(1)$ and $\hat{H}_n=diag\{\hat{h}_n^j\}_{j=0}^p$.

Lemma \ref{lem:norm} simplifies condition (\ref{eq:obj}) further. It allows us to use a nonstochastic normalization in order to obtain the asymptotic properties of the local polynomial estimator constructed with stochastic bandwidths. Throughout the paper, for an arbitrary stochastic vector $W_n$, all orders in probability are taken element-wise.

\begin{lemma}\label{lem:norm}
Define $\Delta_n(h)=Hb_n(x;h)-E(Hb_n(x;h)|\vec{X})$. Suppose that $\sqrt{nh_n}\left(\Delta_n(h_n)-\Delta_n(\hat{h}_n)\right)=o_p(1)$ and $\sqrt{nh_n}\Delta_n(h_n)\stackrel{d}{\rightarrow}W$
a suitably defined random variable. 
Then it follows that
$\sqrt{nh_n}\Delta_n(h_n)-\sqrt{n\hat{h}_n}\Delta_n(\hat{h}_n)=o_p(1)$
provided that $\frac{\hat{h}_n}{h_n}-1=o_p(1)$.
\end{lemma}

Our subsequent results depend on the following assumptions.

\noindent \bf A1. \rm 1. The process $\{(Y_t,X_t)\}_{t=1}^n$ is strictly stationary. 2. for some $\delta>2$ and $a>1-\frac{2}{\delta}$ we assume that
$\sum_{l=1}^\infty l^a\alpha(l)^{1-\frac{2}{\delta}}<\infty$, where $\alpha(l)$ is a mixing coefficient which is defined below.  3. $\sigma^2(x)\equiv Var(Y_t|X_t=x)$ is a continuous and differentiable function at $x$. 4. The $p^{th}$-order derivative of the regressions, $m^{(p)}(x)$, exists at $x$.

The mixing coefficient $\alpha(j)$ is defined as $\alpha(j)\equiv\sup_{A\in\mathcal{F}_{-\infty}^t, B\in\mathcal{F}_{t+j}^\infty}|P(A\cap B)-P(A)P(B)|$ where for a sequence of strictly stationary random vectors $\{(Y_t,X_t)\}_{t\in\mathbb{Z}}$ defined on the probability space $(\Lambda,\mathcal{A},P)$ we define $\mathcal{F}_{a}^b$ as the $\sigma$-algebra induced by $((Y_a,X_a),\ldots,(Y_b,X_b))$ for $a\leq b$ (\cite{Doukhan1994}).  If $\alpha(j)=O(j^{-a-\epsilon})$ for $a \in \mathbb{R}$ and some $\epsilon>0$, $\alpha$ is said to be of size $-a$.  Condition A1.2 is satisfied by a large class of stochastic processes.  In particular, if $\{(Y_t,X_t)\}_{t=1,2,\cdots}$ is $\alpha$-mixing of size $-2$, i.e., $\alpha(l)=O(l^{-2-\epsilon})$ for some $\epsilon>0$ then A1.2 is satisfied.  Since \cite{Pham1985} have shown that finite dimensional stable vector ARMA process are $\alpha$-mixing with $\alpha(l) \rightarrow 0 $ exponentially as $l \rightarrow \infty$, we have that these ARMA processes have size $-a$ for all $a \in \Re^+$, therefore satisfying A1.2.\footnote{Linear stochastic processes also satisfy A1.2 under suitable restrictions.  See \cite{Pham1985} Theorem 2.1.}

\noindent \bf A2. \rm 1. The bandwidth $0<h_n\rightarrow0$ and $nh_n\rightarrow\infty$ as $n\rightarrow\infty$. 2. There exists a stochastic bandwidth  $\hat{h}_n$ such that $\frac{\hat{h}_n}{h_n}-1=o_p(1)$ holds.

\noindent \bf A3. \rm 1. The kernel function $K:\mathbb{R}\rightarrow\mathbb{R}$ is a bounded density function with support $supp(K)=[-1,1]$. 2. $u^{2p\delta+2}K(u)\rightarrow0$ as $|u|\rightarrow\infty$ for $\delta>2$. 3. The first derivative of the kernel function, $K^{(1)}$, exists almost everywhere with $K^{(1)}$uniformly bounded whenever it exists.

\noindent \bf A4. \rm The density $f_X(x)$ for $X_t$ is differentiable and satisfies a Lipschitz condition of order 1, i.e.,
$|f_X(x)-f_X(x')|\leq C|x-x'|$, $\forall x,x'\in\mathbb{R} $.

\noindent \bf A5. \rm 1. The joint density of $(X_t,X_{t+s})$, $f_s(u,v)$, is such that $f_s(u,v)\leq C$ for all $s\geq1$ and $u,v\in[x-h_n,x+h_n]$. 2. $|f_s(u,v)-f_X(u)f_X(v)|\leq C$ for all $s\geq 1$.

\noindent \bf A6. \rm $E(Y_1^2+Y_l^2|X_1=u,X_l=v)<\infty,~~\forall l\geq 1$ and $E(|Y_t|^\delta|X_t=u)<\infty,~~\forall t$, for all $u,v\in[x-h_n,x+h_n]$ and some $\delta>2$.

\noindent \bf A7. \rm There exists a sequence of natural numbers satisfying $s_n\rightarrow\infty$ as $n\rightarrow\infty$ such that
$s_n=o(\sqrt{nh_n})$ and $\alpha(s_n)=o\left(\sqrt{\frac{h_n}{n}}\right)$.

Assumption A7 places a restriction on the speed at which the mixing coefficient decays to zero relative to $h_n$.  Specifically, since the distributional convergence in equation (5) below is established using the large-block/small-block method in \cite{Bernstein1927}, the speed at which the small-block size evolves as $n\rightarrow \infty$ is related to speed of decay for $\alpha$.  In fact, as observed by Masry and Fan (1997), if $h_n \sim n^{-1/5}$ and $s_n=(nh_n)^{1/2}/log n$ it suffices to have $n^g\alpha(n)=O(1)$ for $g>3$ to satisfy A7 (and A1.2).

\noindent \bf A8. \rm The conditional distribution of $Y$ given $X=u$, $f_{Y|X=u}(y)$ is continuous at the point $u=x$.

Let
\begin{eqnarray}\label{sndef}
\label{sndef} s_{n,l}(x;h_n)&=& \frac{1}{nh_n}\sum_{t=1}^nK\left(\frac{X_t-x}{h_n} \right)\left(\frac{X_t-x}{h_n} \right)^l,  
\\
\label{gndef} g_{n,l}(x;h_n)&=& \frac{1}{nh_n}\sum_{t=1}^nK\left(\frac{X_t-x}{h_n} \right)\left(\frac{X_t-x}{h_n} \right)^lY_t ~~~~~~~~~~~ \mbox{and}
\\
\label{gsndef} g_{n,l}^*(x;h_n)&=& \frac{1}{nh_n}\sum_{t=1}^nK\left(\frac{X_t-x}{h_n} \right)\left(\frac{X_t-x}{h_n} \right)^l(Y_t-m(X_t))~~~~~\mbox{for}~l=1,\ldots,2p.
\end{eqnarray}
Then $b_n(x;h_n)=H_n^{-1}S_n^{-1}(x;h_n)G_n(x;h_n)$
where $S_n(x;h_n)=\{s_{n,i+j-2}(x;h_n)\}_{i,j=1}^{p+1,p+1}$ and $G_n(x;h_n)=\{g_{n,l}(x;h_n)\}_{l=0}^p$. \cite{Masry1997} show that under assumption A1 through A8
\begin{equation}\label{mfres}
\sqrt{nh_n}(H_nb_n(x;h_n)-E(S_n(x;h_n)^{-1}G_n(x;h_n)|\vec{X}))\stackrel{d}{\rightarrow}N\left(0,~\frac{\sigma^2(x)}{f_X(x)}S^{-1}\tilde{S}S^{-1}\right),
\end{equation}
where $S=\{\mu_{i+j-2}\}_{i,j=1}^{p+1,p+1}$, $\tilde{S}=\{\nu_{i+j-2}\}_{i,j=1}^{p+1,p+1}$ with $\mu_l=\int\psi^lK(\psi)d\psi$ and $\nu_l=\int\psi^lK^2(\psi)d\psi$.

Equation (\ref{mfres}) gives us  $\sqrt{nh_n}\Delta_n(h_n)\stackrel{d}{\rightarrow}W$ in Lemma \ref{lem:norm}. In particular, $W\sim N(0,~f_X^{-1}(x)\sigma^2(x)S^{-1}\tilde{S}S^{-1})$. Consequently it suffices to show that 
\begin{equation}\label{rdobj}
\sqrt{nh_n}\Delta_n(h_n)-\sqrt{nh_n}\Delta_n(\hat{h}_n)=o_p(1).
\end{equation}

As will be seen in Theorem \ref{thm:main}, the key to establish (\ref{rdobj}) resides in obtaining  asymptotic  uniform stochastic equicontinuity of $\sqrt{nh_n}\Delta_n(x;\tau h_n)$ with respect to $\tau$. To this end we establish the following auxiliary lemmas.
 
\begin{lemma}\label{lem:bdvar1}
Let $Z_n(x;l,\tau)=\left|\frac{d}{d \tau}s_{n,l}(x;\tau h_n)\right|$, for some $\tau$ finite and $l=0,\ldots,2p$. If A1.1, A2.1, A3.1, A3.3 and A4 hold, then $\sup_{\tau\in[r,s]}Z_n(x;l,\tau)=O_p(1)$ where $r,s>0$ and $r<s$.
\end{lemma}
\begin{lemma}\label{lem:bdvar2}
Let $B_n(x;l,\tau)=\sqrt{nh_n}\frac{d}{d\tau}g_{n,l}^*(x;\tau h_n)$, for $l=0,\ldots,2p$. If A1 through A6 hold, then $\int_r^sB_n^2(x,l,\tau)d\tau=O_p(1)$
where $r,s>0$ and $r<s$.
\end{lemma}

\section{Main Results}
\label{sec:lpmr}

The following theorem and corollary establish $\sqrt{n\hat{h}_n}$-normality of the local polynomial estimator constructed with stochastic bandwidths. As in \cite{Masry1997}, we are able to obtain asymptotic normality for the regression estimator as well as for the estimators of the regression derivatives.
\begin{theorem}\label{thm:main}
Suppose A1 through A8 hold, then it follows that
 $$\sqrt{n\hat{h}_n}\Delta_n(\hat{h}_n)\stackrel{d}{\rightarrow}N\left(0,\frac{\sigma^2(x)}{f_X(x)}S^{-1}\tilde{S}S^{-1}\right).$$
\end{theorem}

With the following corollary we also obtain the asymptotic bias for local polynomial estimators with stochastic bandwidths.

\begin{corollary}\label{col:main}
Let $m^{(j)}$ denote the $j$th-order derivative of $m$. Suppose A1 through A8 hold, then
$$\sqrt{n\hat{h}_n}\left(\hat{H}_n(b_n(x;\hat{h}_n)-b(x))-\frac{\hat{h}_n^{p+1}m^{(p+1)}(x)}{(p+1)!}+\hat{h}_n^{p+1}o_p(1) \right)\stackrel{d}{\rightarrow}N\left(0,\frac{\sigma^2(x)}{f_X(x)}S^{-1}\tilde{S}S^{-1}\right)$$
where $\hat{H}_n=diag\{\hat{h}_n^j\}_{j=0}^p$ and $b(x)=\left(m(x), m^{(1)}(x),\ldots,\frac{1}{p!}m^{(p)}(x) \right)'$.
\end{corollary}

\section{Monte Carlo study}
In this section we investigate some of the finite sample properties of the local linear regression and derivative estimators constructed with bandwidths selected by cross validation and a plug in method proposed by \cite{Ruppert1995} for data generating processes (DGP) exhibiting dependence. In our simulations two regression functions are considered, $m_1(x)=\sin(x)$ and $m_2(x)=3(x-0.5)^3+0.25x_0+1.125$ with first derivatives given respectively by $m_1^{(1)}(x)=-\cos(x)$ and $m_2^{(1)}(x)=9(x-0.5)^2+0.25$.

We generate $\{\epsilon_t\}_{t=1}^n$ by $\epsilon_t=\rho\epsilon_{t-1}+\sigma U_t$, where $\{U_t\}_{t\geq1}$ is a sequence of IID standard normal random variable and $(\rho, \sigma^2)=(0, 0.01),~(0.5, 0.0075)$ and $(0.9, 0.0019)$. This implies that for $\rho \neq 0$ $\{\epsilon_t\}_{t=1}^n$ is a normally distributed AR(1) process with mean zero and variance equal to 0.01.

For $m_1$ we draw IID regressors $\{X_t\}_{t=1}^n$ from a uniform distribution that takes value on $[0,2\pi]$.
For $m_2$ we draw IID regressors $\{X_t\}_{t=1}^n$ from a beta distribution with parameters $\alpha=2$ and $\beta=2$ given by 
\begin{equation}\notag
f_X(x;\alpha,\beta)=\left\{
\begin{array}{ll}
\frac{x^{\alpha-1}(1-x)^{\beta-1}}{\int_0^1u^{\alpha-1}(1-u)^{\beta-1}du} & \mbox{if}~x\in[0,1] \\
0 & \mbox{otherwise.}
\end{array}
\right. 
\end{equation} 
 The regressands are constructed using $Y_t=m_i(X_t)+\epsilon_t$, where $i=1,2$.

Two sample sizes are considered $n=200, 600$ and $1000$ repetitions are performed. We evaluate the regression and regression derivative estimators at $x=0.5\pi,~\pi,~1.5\pi$ and $x=0.25,~0.5,~0.75$ for $m_1$ and $m_2$ respectively.  These estimators are constructed with a nonstochastic optimal regression bandwidth $h_{AMISE}$ and with two data dependent bandwidths: a cross validated bandwidth $h_{CV}$ and a plug in bandwidth $h_{ROT}$.
The nonstochastic bandwidth is given by $h_{AMISE}=\left(\frac{\lambda_1}{n\lambda_2} \right)^{1/5}$, where $\lambda_1=Var(\epsilon_t)\int K^2(u)du\int \mathbf{1}(f_X(x)\neq 0)dx$ and $\lambda_2=\int u^2K(u)du\int(m^{(2)}(x))^2f_X(x)dx$ (\cite{Ruppert1995}).  The cross validated bandwidth is given by $h_{CV}=\underset{h}{\mbox{argmin}}\sum_{t=1}^n(m_{LP,t}(X_t;h)-Y_t)^2$, where $m_{LP,t}(x;h)$ is the local linear regression estimator constructed with the exclusion of observation $t$ (\cite{Xia2002}).  The $h_{ROT}$ bandwidth is calculated as described in \cite{Ruppert1995}.  Specifically, we estimate 
$Var(\epsilon_t)$, $\int \mathbf{1}(f_X(x)\neq 0)dx$ and $\int(m^{(2)}(x))^2f_X(x)dx$ which appear the expression for $h_{AMISE}$.  First, we approximate $m(x)$ by $m(x)\approx \beta_0+\beta_1x+\beta_2x^2/2+\beta_3x^3/3!+\beta_4x^4/4!$ and obtain $m^{(2)}(x)\approx\beta_2+\beta_3x+\beta_4x^2/2$.  Second, the vector $(\beta_0,\ldots,\beta_4)'$ is estimated by $(\hat{\beta}_0,\ldots,\hat{\beta}_4)'=(\sum_{i=1}^nR_{i\cdot}'R_{i\cdot})^{-1}\sum_{i=1}^nR_{i\cdot}'Y_i$, where $R_{i\cdot}=(1~X_i~X_i^2/2 \cdots X_i^4/4!)$.  Third, we estimate $Var(\epsilon_t)$ and $\int(m^{(2)}(x))^2f_X(x)dx$ by $n^{-1}\sum_{i=1}^n\tilde{e}_i^2$ and $n^{-1}\sum_{i=1}^n[\hat{\beta}_2+\hat{\beta}_3X_i+\hat{\beta}_4X_i^4/2]^2$ respectively, where $\tilde{e}_i=y_i-R_{i\cdot}(\hat{\beta}_0,\ldots,\hat{\beta}_4)'$. The estimator used for $\int \mathbf{1}(f_X(x)\neq 0)dx$ is given by $\max_iX_i-\min_iX_i$.

\begin{figure}[h]
\caption{Estimated density of regression}\label{fig:denreg}
\centerline{\begin{tabular}{cc}
(a) Estimated density of $m_1(0.5\pi)$ using $h_{CV}$ & (b) Estimated density of $m_1(0.5\pi)$ using $h_{ROT}$ \\
\includegraphics[width=3.5in]{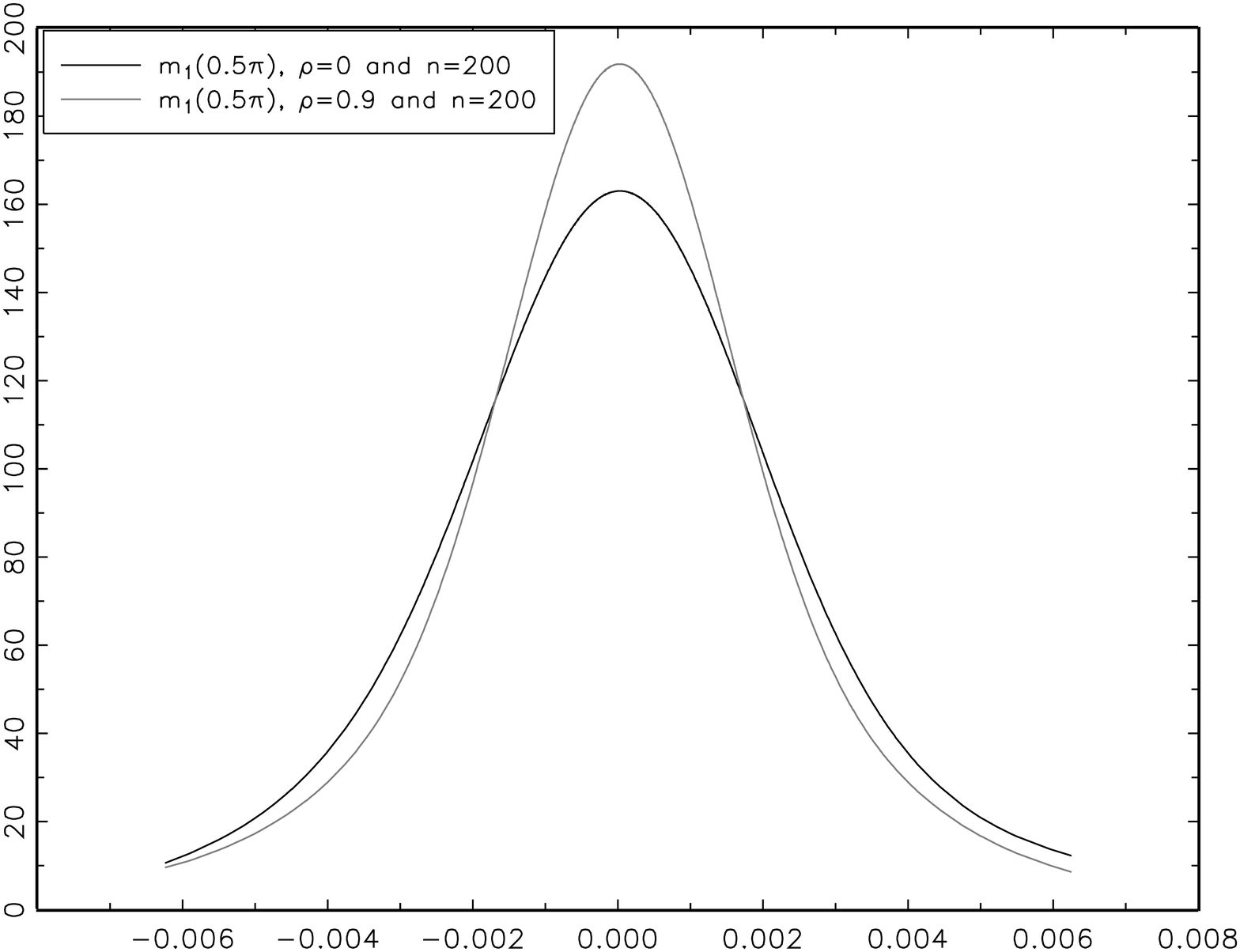} & \includegraphics[width=3.5in]{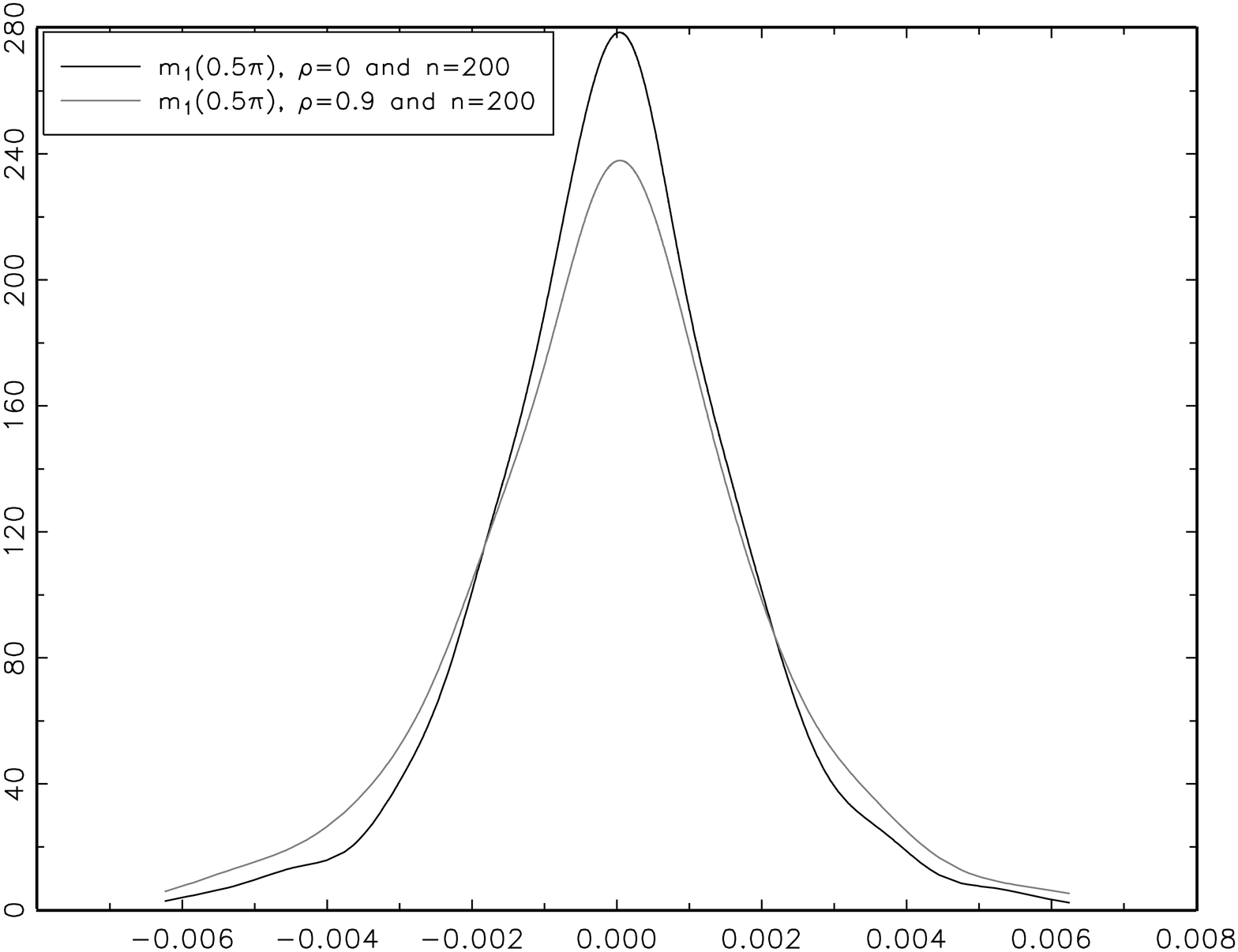}
\\ \\
(c) Estimated density of $m_2(0.5)$ using $h_{CV}$ & (d) Estimated density of $m_2(0.5)$ using $h_{ROT}$\\
\includegraphics[width=3.5in]{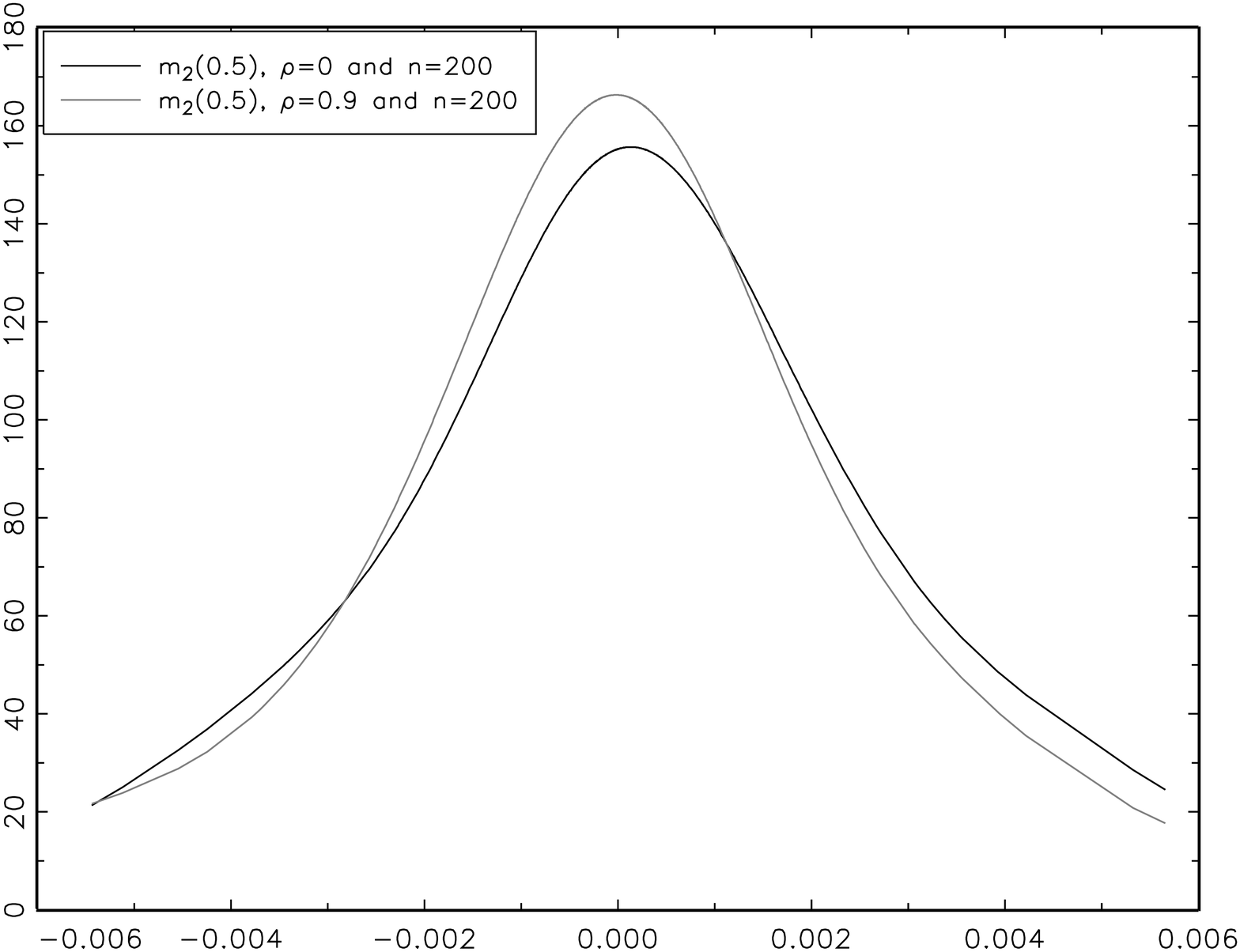} & \includegraphics[width=3.5in]{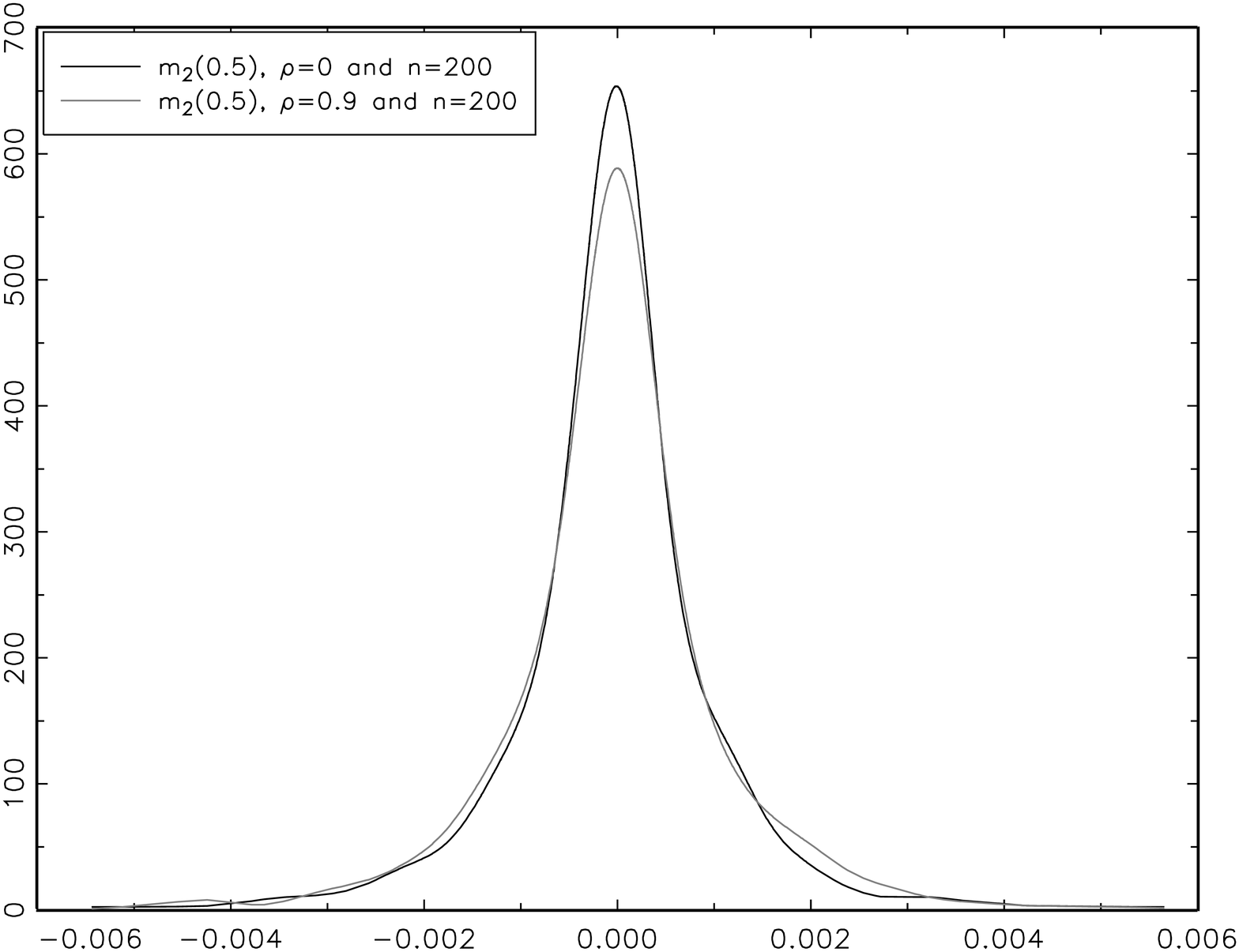}
\end{tabular}
}
\end{figure}

\begin{figure}[h]
\caption{Estimated density of regression derivative}\label{fig:dender}
\centerline{\begin{tabular}{cc}
(a) Estimated density of $m_1^{(1)}(0.5\pi)$ using $h_{CV}$ & (b) Estimated density of $m_1^{(1)}(0.5\pi)$ using $h_{ROT}$\\
\includegraphics[width=3.5in]{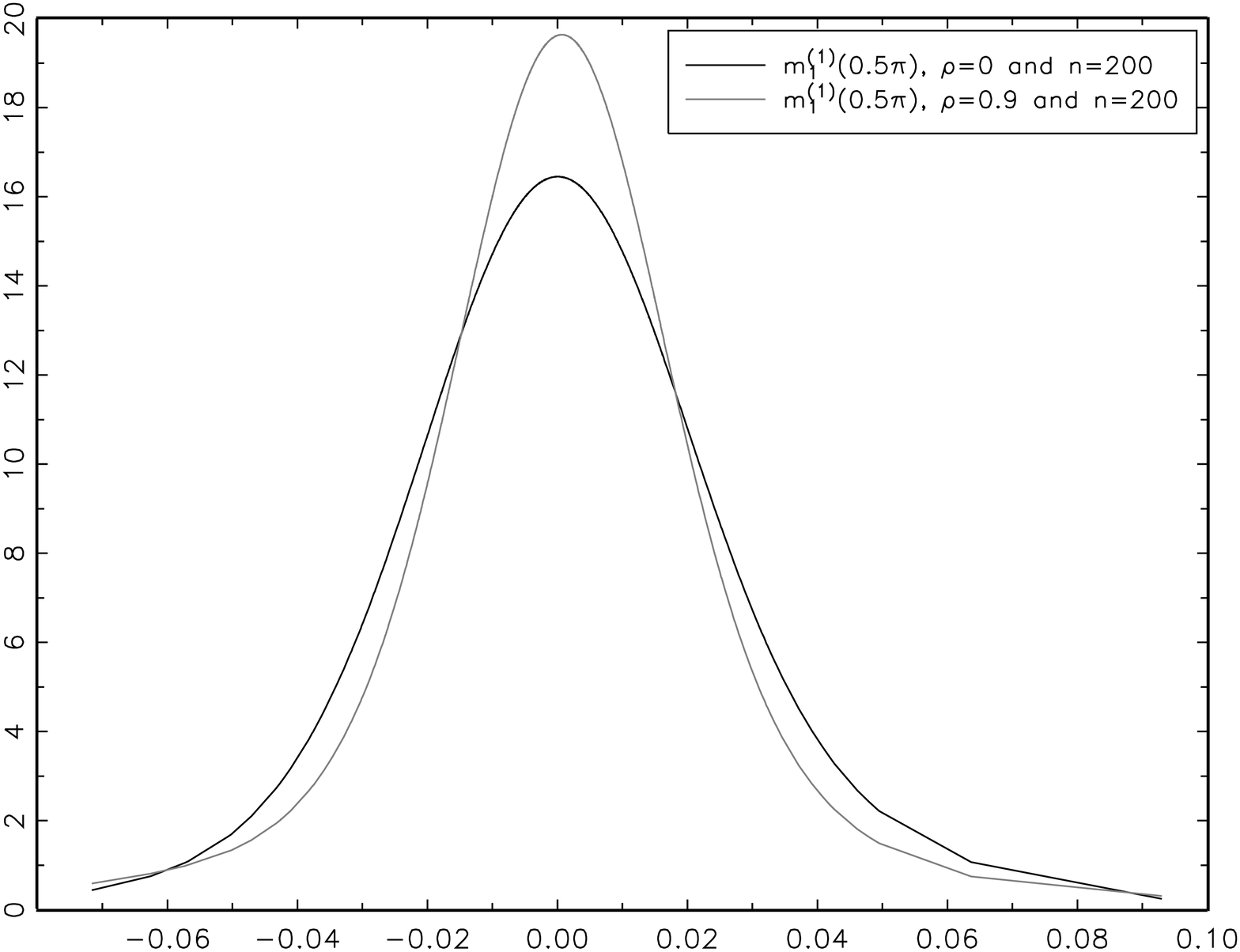} & \includegraphics[width=3.5in]{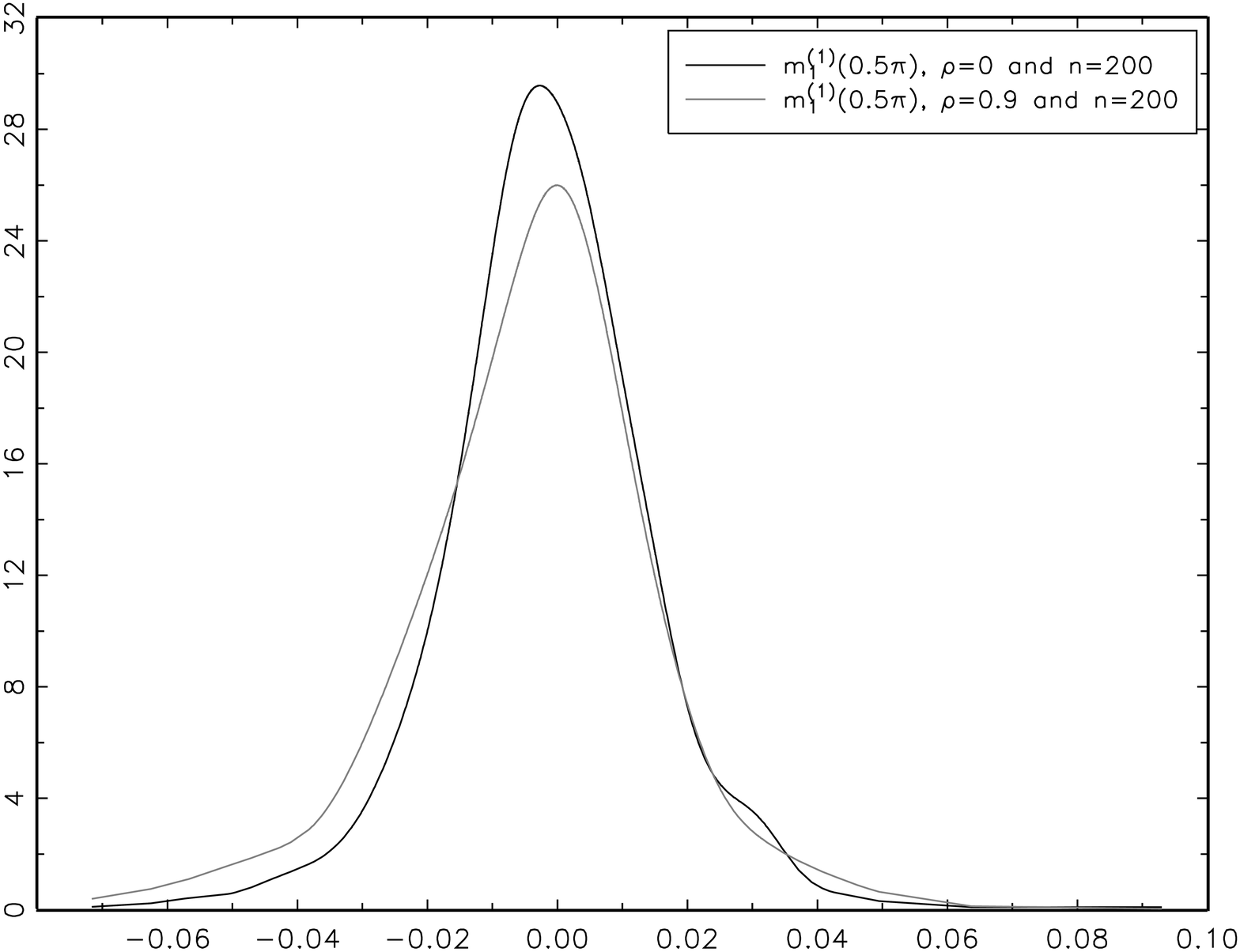} \\
\\
(c) Estimated density of $m_2^{(1)}(0.5)$ using $h_{CV}$ & (d) Estimated density of $m_2^{(1)}(0.5)$ using $h_{ROT}$\\
\includegraphics[width=3.5in]{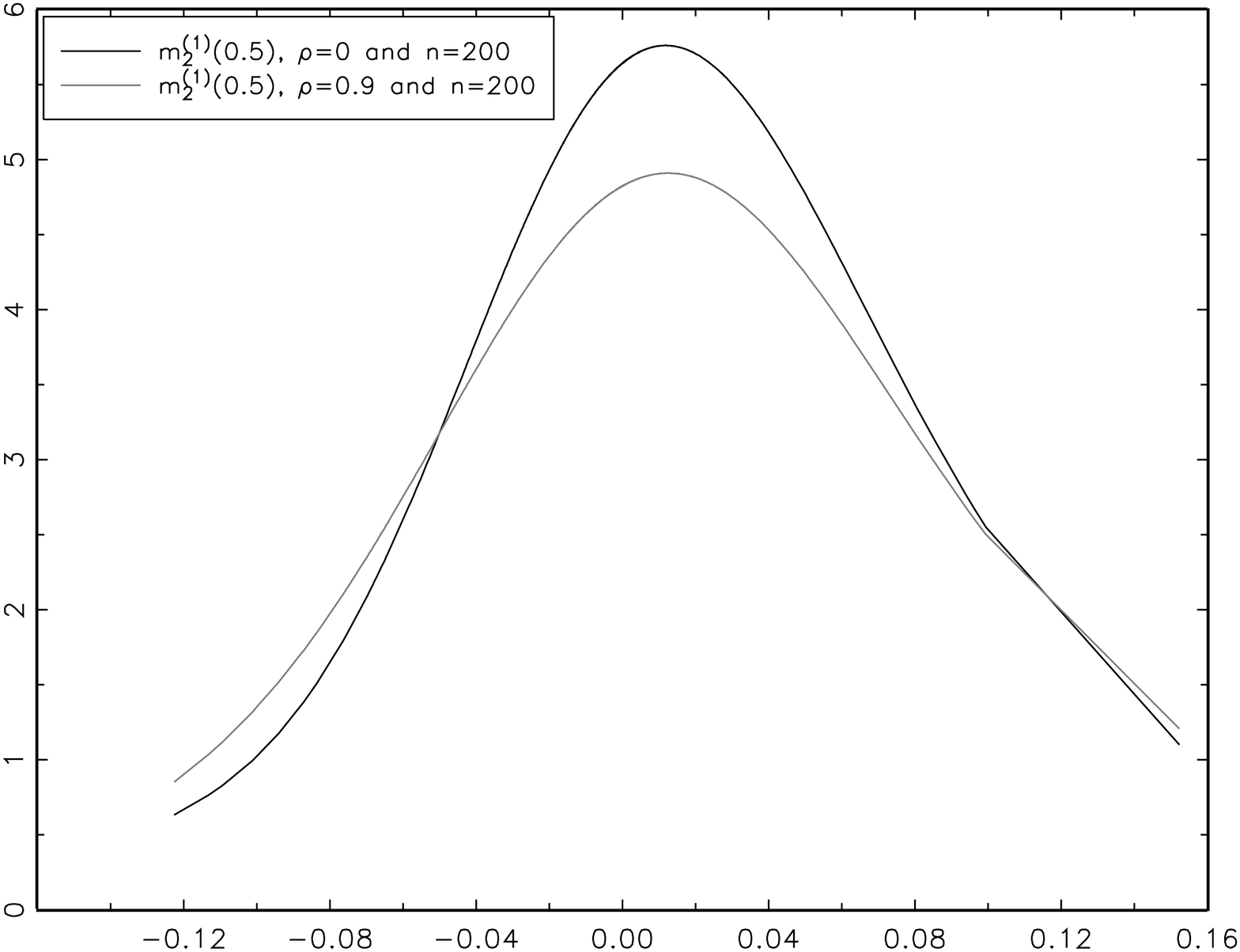} & \includegraphics[width=3.5in]{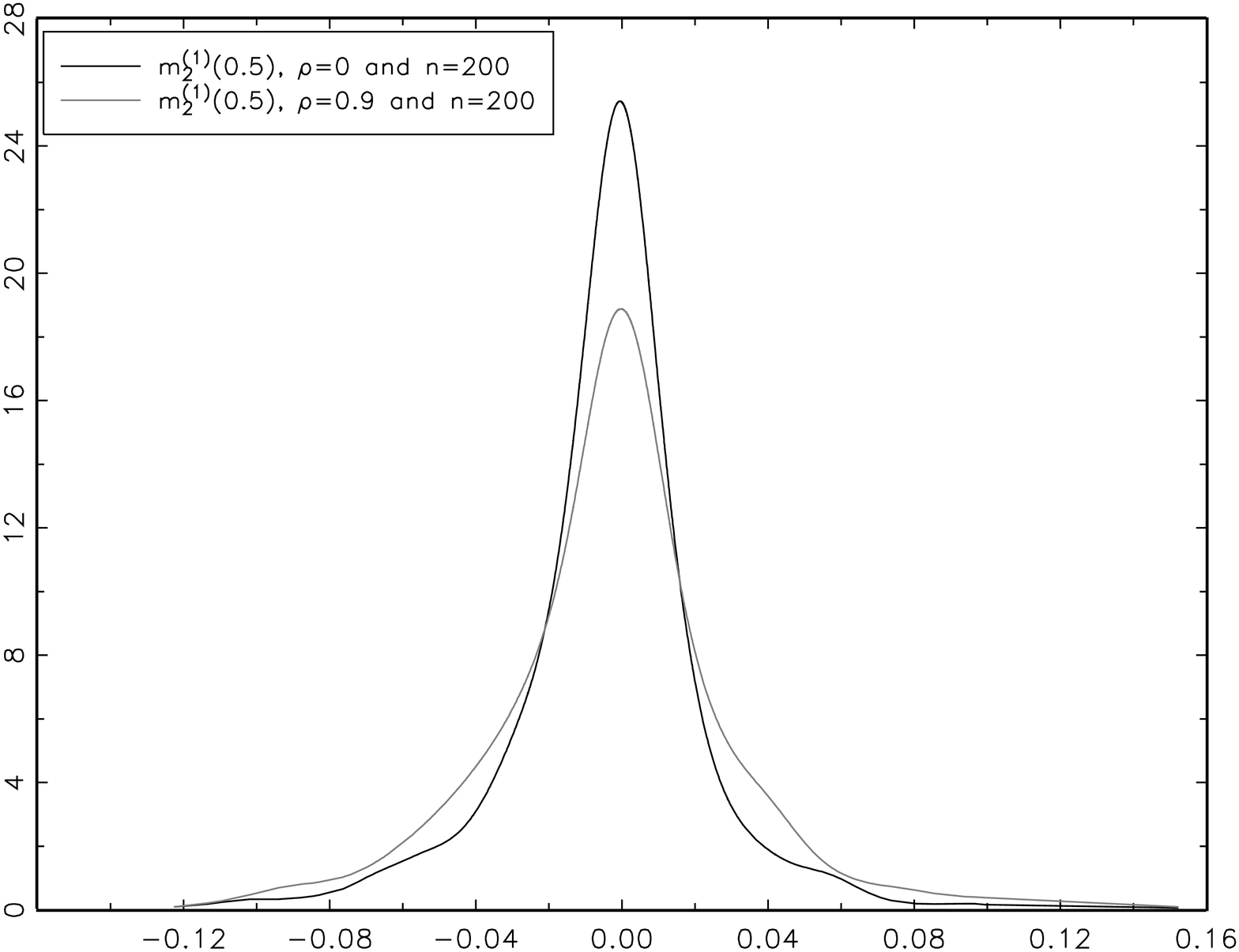}
\end{tabular}
}
\end{figure}

The results of our simulations are summarized in Tables 1-2 and Figures 1-2. Tables \ref{tab:mc1} and \ref{tab:mc2} provide the bias ratio and mean squared error (MSE) ratio of estimators constructed with $h_{CV}$, $h_{ROT}$ and $h_{AMISE}$ for $m_1$ and $m_2$ respectively. These ratios are constructed with estimators using the data dependent bandwidth $h_{CV}$ or $h_{ROT}$ in the numerator and $h_{AMISE}$ in the denominator. Figure \ref{fig:denreg} shows the estimated densities of the difference between the estimated regression constructed with $h_{CV}$ and $h_{AMISE}$ (panels (a) and (c)) and $h_{ROT}$ and $h_{AMISE}$ (panels (b) and (d)), for $m_1(\pi)$ and $m_2(0.5)$, with $n=200$ and $\rho=0,0.9$. Similarly, figure \ref{fig:dender} shows the estimated density of the difference between the estimated regression first derivative constructed with $h_{CV}$ and $h_{AMISE}$ (panels (a) and (c)) and $h_{ROT}$ and $h_{AMISE}$ (panels (b) and (d)) for $m_1^{(1)}(\pi)$ and $m_2^{(1)}(0,5)$, with $n=200$ and $\rho=0,0.9$.    

As expected from the asymptotic results, the bias and MSE ratios are in general close to 1, especially for the regression estimators.  Ratios that are farther form 1 are more common in the estimation of the regression derivatives. This is consistent with the asymptotic results since the rate of convergence of the regression estimator is $\sqrt{nh_n}$, whereas regression first derivative estimators have rate of convergence $\sqrt{nh_n^3}$.  Hence, for fixed sample sizes we expect regression estimators to outperform those associated with derivatives.

Note that most bias and MSE ratios given in tables \ref{tab:mc1} and \ref{tab:mc2} are positive values larger than 1. Since we constructed both bias and MSE ratios with estimators constructed with $h_{CV}$ or $h_{ROT}$ in the numerator and estimators constructed with $h_{AMISE}$ in the denominator, the results indicate that bias and MSE are larger for estimators constructed with $h_{CV}$ and $h_{ROT}$. This too was expected, since $h_{AMISE}$ is the true optimal bandwidth for the regression estimator. Positive bias ratios indicate that the direction of the bias is the same for estimators constructed with $h_{AMISE}$ and $h_{CV}$ or $h_{ROT}$.  It is also important to note that the bias and MSE for estimators of both regression and derivatives are generally larger when calculated using $h_{CV}$ compared to the case when $h_{ROT}$ is used.   

We note that  in general the estimators for the function $m_1$ outperformed those for function $m_2$.  We observe that $m_2$ takes value on $[0.75,1.75]$ and $\epsilon_t$ on $\mathbb{R}$. Thus, although the variance of $\epsilon_t$ was chosen to be small, $0.01$, estimating the bandwidth was made difficult due to the fact that $\epsilon_t$ had a large impact on $Y_t$ in terms of its relative magnitude. The regression function $m_1$ also took values on a bounded interval, however this interval had a larger range. In fact the standard deviation of $h_{CV}$ for $n=200$ and $\rho=0.5$ was $0.0422$ and $11.552$ for the DGP's associated with $m_1$ and $m_2$ respectively.

The kernel density estimates shown in figures \ref{fig:denreg} and \ref{fig:dender} were calculated using the Gaussian kernel and bandwidths were selected using the \emph{rule-of-thumb} procedure of \cite{Silverman1986}. We observe that the change from IID ($\rho=0$) to dependent DGP ($\rho\neq 0$) did not yield significantly different results in terms of estimator performance under $h_{CV}$ or $h_{ROT}$.  In fact, our results seem to indicate that for $\rho=0.9$ the estimators had slightly better general performance than for the case where $\rho=0$.  As expected from our asymptotic results, figures \ref{fig:denreg} and \ref{fig:dender} show that the difference between derivative estimates using $h_{CV}$ and $h_{AMISE}$ and $h_{ROT}$ and $h_{AMISE}$ were more disperse around zero than those associated with regression estimates, especially for the DGP using $m_2$. Even though the DGP for $m_1$ provided better results, the estimators of $m_2$ and $m_2^{(1)}$, as seen on figures \ref{fig:denreg} panels (c) and (d) and \ref{fig:dender} panels (c) and (d) performed well, in the sense that such estimators produced estimated densities with fairly small dispersion around zero.  Another noticeable result from the Monte Carlo is that the estimated densities associated with estimators calculated using $h_{ROT}$ are much less dispersed than those calculated using $h_{CV}$.  Overall, as expected from asymptotic theory, estimators calculated with $h_{CV}$ and $h_{ROT}$ performed fairly well in small samples, however our results seem to indicate better performance when a plug in bandwidth is used.

\section{Final Remarks}

We have established the asymptotic properties of the local polynomial regression estimator constructed with stochastic bandwidths. Our results validate the use of the normal distribution in the implementation of hypotheses tests and interval estimation when bandwidths are data dependent. Most assumptions that we have imposed, were also explored by \cite{Masry1997}. The assumptions we place on $\hat{h}_n$ coincides with the properties of the bandwidths proposed by \cite{Ruppert1995} and \cite{Xia2002} under IID and strong mixing respectively.

\section*{Appendix 1: Proofs}

\textbf{Proof of Lemma \ref{lem:norm}:}
\begin{proof}
Since $\sqrt{nh_n}(\Delta_n(h_n)-\Delta_n(\hat{h}_n))=o_p(1)$, we have that
$$\sqrt{nh_n}\Delta_n(h_n)-\sqrt{n\hat{h}_n}\Delta_n(\hat{h}_n)=o_p(1)+\sqrt{nh_n}\Delta_n(\hat{h}_n)-\sqrt{n\hat{h}_n}\Delta_n(\hat{h}_n).$$
$\sqrt{nh_n}\Delta_n(h_n)\stackrel{d}{\rightarrow}W$ and $\sqrt{nh_n}\Delta_n(h_n)-\sqrt{nh_n}\Delta_n(\hat{h}_n)=o_p(1)$ imply that $\Delta_n(\hat{h}_n)=O_p\left((nh_n)^{-1/2}\right)$. Consequently,
$$\sqrt{nh_n}\Delta_n(\hat{h}_n)-\sqrt{n\hat{h}_n}\Delta_n(\hat{h}_n)=\left(1-\sqrt{\frac{\hat{h}_n}{h_n}}\right)O_p(1)=o_p(1)$$
since $\left(1-\sqrt{\frac{\hat{h}_n}{h_n}}\right)=o_p(1)$. 
\end{proof}

\noindent\textbf{Proof of Lemma \ref{lem:bdvar1}:}
\begin{proof}
For any $\epsilon>0$, we must find $M_\epsilon<\infty$ such that
\begin{equation}
P\left(\sup_{\tau\in [r,s]}Z_n(x;l,\tau)>M_\epsilon\right)\leq\epsilon.
\end{equation}
By Markov's inequality, we have that
\begin{equation}
P\left(\sup_{\tau\in [r,s]}Z_n(x;l,\tau)>\frac{1}{\epsilon}\right)\leq E\left(\sup_{\tau\in [r,s]}Z_n(x; l,\tau)\right)\epsilon.
\end{equation}
Thus it suffices to show that $E\left(\sup_{\tau\in [r,s]}Z_n(x;l,\tau)\right)=O(1)$. Let $\tilde{K}_l(x)=K(x)x^l(1+l)+K^{(1)}(x)x^{l+1}$ and write
\begin{equation}
Z_n(x;l,\tau)=\frac{1}{nh_n\tau^2}\left|\sum_{t=1}^n\tilde{K}_l\left(\frac{X_t-x}{h_n\tau}\right)\right|.
\end{equation}
By strict stationarity, we write
\begin{equation}
E\left(\sup_{\tau\in [r,s]}Z_n(x;l,\tau)\right)\leq\frac{1}{h_nr^2}E\left(\sup_{\tau\in [r,s]}\left|\tilde{K}_l\left(\frac{X_t-x}{h_n\tau}\right)\right|\right)
\end{equation}
Now, note that
\begin{equation}
\begin{array}{rcl}
E\left(\sup_{\tau\in [r,s]}\left|\tilde{K}_l\left(\frac{X_t-x}{h_n\tau}\right)\right|\right)&=&h_n\int\sup_{\tau\in [r,s]}|\tau\tilde{K}_l(\phi)||f_X(x+h_n\tau\phi)-f_X(x)+f_X(x)|d\phi\\
&\leq& h_n^2C\int|\tilde{K}_l(\phi)\phi|d\phi\sup_{\tau\in [r,s]}\tau^2+h_nf_X(x)\int|\tilde{K}_l(\phi)|d\phi\sup_{\tau\in [r,s]}\tau\\
&\leq& h_n^2s^2C\int|\tilde{K}_l(\phi)\phi|d\phi\sup_{\tau\in [r,s]}+h_nf_X(x)s\int|\tilde{K}_l(\phi)|d\phi\\
&\leq&h_n^2s^2C\int_{-1}^1(|1+l||K(\phi)||\phi^{l+1}|+|K^{(1)}(\phi)||\phi^{l+2}|)d\phi\\
& &+h_nf_X(x)s\int_{-1}^1(|1+l||K(\phi)||\phi^{l}|+|K^{(1)}(\phi)||\phi^{l+1}|)d\phi\\
&\leq &(h^2_ns^2C+sh_nf_X(x))\int_{-1}^1(|1+l||K(\phi)|+|K^{(1)}(\phi)|)d\phi
\end{array}
\end{equation}
Hence,
\begin{eqnarray}
E\left(\sup_{\tau\in[r,s]}Z_n(x;l\tau)\right)&\leq&\frac{1}{r^2}(h_nC+f_X(x)C)\\
&=&O(1)
\end{eqnarray}
as $h_n\rightarrow0$ and $n\rightarrow\infty$.
\end{proof}

\noindent\textbf{Proof of Lemma \ref{lem:bdvar2}:}
\begin{proof}
Using Markov's inequality it suffices to establish that
\begin{equation}
E\left(\int_r^sB_n^2(x;l,\tau)d\tau\right)=\int_r^sE\left(B_n^2(x;l,\tau)\right)d\tau=O(1).
\end{equation}
Note that,
\begin{equation}
B_n^2(x;l,\tau)=\frac{1}{nh_n\tau^4}\left\{\sum_{t=1}^n\tilde{K}_l^2\left(\frac{X_t-x}{h_n\tau}\right)\epsilon_t^2+2\sum_{t=1}^n\sum_{i\neq t}\tilde{K}_l\left(\frac{X_t-x}{h_n\tau}\right)\epsilon_t\tilde{K}_l\left(\frac{X_i-x}{h_n\tau}\right)\epsilon_i\right\}
\end{equation}
where $\epsilon_t=Y_t-m(X_t)$. Thus, by the law of iterated expectations and strict stationarity, we obtain
\begin{equation}
\begin{array}{rl}
E\left(B_n^2(x;l,\tau)\right)= &\left|E\left(\frac{1}{h_n\tau^4}\tilde{K}_l^2\left(\frac{X_t-x}{h_n\tau}\right)\sigma^2(X_t)\right)\right.\\
+&\left.2\frac{1}{h_n\tau^4}\sum_{t=2}^n\left(1-\frac{t}{n}\right)E\left(\tilde{K}_l\left(\frac{X_1-x}{h_n\tau}\right)\tilde{K}_l\left(\frac{X_t-x}{h_n\tau}\right)\epsilon_1\epsilon_t\right)\right|\\
\leq&E\left(\frac{1}{h_n\tau^4}\tilde{K}_l^2\left(\frac{X_t-x}{h_n\tau}\right)\sigma^2(X_t)\right)\\
+&2\frac{1}{h_n\tau^4}\sum_{t=2}^n\left(1-\frac{t}{n}\right)\left|E\left(\tilde{K}_l\left(\frac{X_1-x}{h_n\tau}\right)\tilde{K}_l\left(\frac{X_t-x}{h_n\tau}\right)\epsilon_1\epsilon_t\right)\right|
\end{array}
\end{equation}

Notice that 
\begin{equation}
\begin{array}{rcl}
E\left(\frac{1}{h_n\tau^3}\tilde{K}_l^2\left(\frac{X_t-x}{h_n\tau}\right)\sigma^2(X_t)\right)&=&\int\frac{1}{\tau^3}\tilde{K}_l^2(\phi)\sigma^2(x+h_n\tau\phi)f_X(x+h_n\tau\phi)d\phi\\
	&=&\int\tau^{-3}\tilde{K}_l^2(\phi)\left\{\sigma^2(x)f_X(x)+\frac{dw(x^*)}{d x}h_n\tau\phi\right\}d\phi\\
	&\leq&\sigma^2(x)f_X(x)\tau^{-3}\int\tilde{K}_l^2(\phi)d\phi+\tau^{-2}O(h_n)=O(1)
\end{array}
\end{equation}
where $w(x)=f_X(x)\sigma^2(x)$. Let $\xi_t=\tilde{K}_l\left(\frac{X_t-x}{h_n\tau}\right)$, and without loss of generality take $s\geq1$.\footnote{Let $s*=\max\{1,s\}$ and note that this proof follows with $s=s^*$.} Then,
\begin{equation}
\begin{array}{rcl}
|E(\xi_1\xi_t\epsilon_1\epsilon_t)|&=&\left|E\left(E\left(\epsilon_1\epsilon_t|X_1,X_t\right)\xi_1\xi_t\right)\right|\\
&\leq&E\left(\sup_{X_1,X_t\in[x-sh_n,x+sh_n]}E(|\epsilon_1\epsilon_t||X_1,X_t)|\xi_1\xi_t|\right)\\
&\leq&E\left(\sup_{X_1,X_t\in[x-sh_n,x+sh_n]}E\left((|Y_1|+B)(|Y_t|+B)|X_1,X_t\right)|\xi_1\xi_t|\right)\\
&\leq&E\left(\sup_{X_1,X_t\in[x-sh_n,x+sh_n]}\left\{E\left((|Y_1|+B)^2|X_1,X_t\right)E\left((|Y_t|+B)^2|X_1,X_t\right)\right\}^{\frac{1}{2}}|\xi_1\xi_t|\right)\\
&\leq& CE(|\xi_1\xi_t|)\\
&=&C \int\int \left|\tilde{K}_j\left(\frac{u-x}{\tau h_n}\right)\tilde{K}_j\left(\frac{v-x}{\tau h_n}\right)\right|f_t(u,v)dudv\\
&\leq&C \int \int \left|\tilde{K}_j\left(\frac{u-x}{\tau h_n}\right)\tilde{K}_j\left(\frac{v-x}{\tau h_n}\right)\right|dudv\\
&\leq& h_n^2\tau^2C\left(\int\left|\tilde{K}_l(\phi)\right|d\phi\right)^2
\end{array}
\end{equation}
where $B=\sup_{X\in[x-sh_n,x+sh_n]}|m(X)|$. Let $\{d_n\}_{n\geq1}$ be a sequence of positive integers, such that $d_n\rightarrow\infty$ as $n\rightarrow\infty$. Then we can write
\begin{equation}
\sum_{t=2}^n\left|E\left(\xi_1\xi_t\epsilon_1\epsilon_t\right)\right|=\sum_{t=2}^{d_n+1}\left|E\left(\xi_1\xi_t\epsilon_1\epsilon_t\right)\right|+\sum_{t=d_n+2}^n \left|E\left(\xi_1\xi_t\epsilon_1\epsilon_t\right)\right|
\end{equation}
and note that
\begin{equation}
\begin{array}{rcl}
\sum_{t=2}^{d_n+1}\left|E\left(\xi_1\xi_t\epsilon_1\epsilon_t\right)\right|&\leq&\sum_{t=2}^{d_n+1}\tau^2h_n^2C\left(\int\left|\tilde{K}_l(\phi) \right| d\phi\right)^2\\
&=&d_nh_n^2\tau^2C.
\end{array}
\end{equation}

Then using the fact that $E(\xi_t\epsilon_t)=0$ and Davydov's Inequality we obtain,
\begin{equation}
\left|E\left(\xi_1\xi_{t+1}\epsilon_1\epsilon_{t+1}\right)\right|\leq8[\alpha(t)]^{1-2/\delta}\left(E|\xi_1\epsilon_1|^\delta\right)^{2/\delta}.
\end{equation}
Note also that
\begin{equation}
\begin{array}{rcl}
E\left|\xi_1\epsilon_1\right|^\delta&=&E\left|(Y_1-m(X_1))\tilde{K}_l\left(\frac{X_1-x}{h_n\tau}\right)\right|^\delta\\
&\leq&E\left(\sup_{X_1\in[x-sh_n,x+sh_n]}E\{(|Y_1|+B)^\delta|X_1\}\left|\tilde{K}_l\left(\frac{X_1-x}{h_n\tau}\right)\right|^\delta\right) \\
&\leq&CE\left(\left|\tilde{K}_l\left(\frac{X_1-x}{h_n\tau}\right)\right|^\delta\right)\\
&=&C\int\left|\tilde{K}_l\left(\frac{u-x}{h_n\tau}\right)\right|^\delta f_X(u)du\\
&=&C\tau h_n \int |\tilde{K}_l(v)|^\delta f_X(v+\tau h_nx)dv\\
&\leq&C\tau h_n
\end{array}
\end{equation}
which leads to,\
\begin{equation}
\begin{array}{rcl}
\sum_{t=d_n+2}^n\left|E(\xi_1\epsilon_1\xi_{t+1}\epsilon_{t+1})\right|&\leq& \sum_{t=d_n+2}^n8\alpha(t)^{1-2/\delta}\left(E|\xi_1\epsilon_1|^\delta\right)^{2/\delta}\\
&\leq&h_n^{2/\delta}\tau^{2/\delta}C\sum_{t=d_n+2}^n\alpha(t)^{1-2/\delta}\\
&\leq&h_n^{2/\delta}\tau^{2/\delta}C\sum_{t=d_n+2}^n\frac{t^a}{d_n^{1-2/\delta}}\alpha(t)^{1-2/\delta}\\
&=&Ch_n^{2/\delta}d_n^{-1+2/\delta}\tau^{2/\delta}\sum_{t=d_n+2}^nt^a\alpha(t)^{1-2/\delta}\\
&=&Ch_n^{2/\delta}d_n^{-1+2/\delta}\tau^{2/\delta} o(1)\\
&=&C\tau^{2/\delta} o(h_n)
\end{array}
\end{equation}
given that $d_n$ is chosen as the integer part of $h_n^{-1}$ and $a>1-\frac{2}{\delta}$.

Consequently,
\begin{equation}
\begin{array}{rl}
E(B_n^2(x;l,\tau))&=\left(\tau^{-3}O(1)+\tau^{-2}O(h_n)\right)+\left(\tau^{-2}O(1)+\tau^{-4+2/\delta}o(1)\right).
\end{array}
\end{equation}
\end{proof}

\noindent\textbf{Proof of Theorem \ref{thm:main}:}
\begin{proof}
From \cite{Masry1997} and Lemma \ref{lem:norm}, it suffices to show that
$$\sqrt{nh_n}\Delta_n(h_n)-\sqrt{nh_n}\Delta_n(\hat{\tau}_nh_n)=o_p(1)$$
where $\hat{\tau}_n=\frac{\hat{h}_n}{h_n}$. It suffices to show that all elements of the vector $\sqrt{nh_n}\Delta_n(\tau h_n)$ are stochastically equicontinuous on $\tau$.

For any $\epsilon>0$, given that $\hat{\tau}_n=O_p(1)$ there exists $r,s\in(0,\infty)$ with $r<s$ such that
$P(\hat{\tau}_n\notin[r,s])\leq\epsilon/3,~~~\forall~n.$
For $\delta>0$, let
$w_n(i,\delta)=\sup_{\{(\tau_1,\tau_2)\in[r,s]\times[r,s]:|\tau_1-\tau_2|<\delta \}}d_n(i,\tau_1,\tau_2)$
where
$d_n(x;i,\tau_1,\tau_2)=|e_i\Delta_n(\tau_2h_n)-e_i\Delta_n(\tau_1h_n)|,$
$e_i$ is a row vector with $i$-th component equal to 1, and 0 elsewhere. Then, for $\eta>0$
$$P(d_n(i,1,\hat{\tau}_n)\geq\eta)\leq P(\mathbf{1}(|\hat{\tau}_n-1|\leq\delta)d_n(i,1,\hat{\tau}_n)\geq\eta)+ P(|\hat{\tau}_n-1|>\delta)
$$
where $1(A)$ is the indicator function for the set $A$.

By assumption, there exists $N_{\epsilon,1}$ such that $P(|\hat{\tau}_n-1|>\delta)\leq\frac{\epsilon}{3}$, $\forall~ n\geq N_{\varepsilon,1}$. Also,

\begin{equation}
\begin{array}{lcl}
P(\mathbf{1}(|\tau_n-1|\leq\delta)d_n(i,1,\hat{\tau}_n)\geq\eta)&\leq &P(\mathbf{1}(\hat{\tau}_n\in[1-\delta,1+\delta]\cap[r,s])d_n(i,1,\hat{\tau}_n)\geq\eta)\\& +&P(\hat{\tau}_n\notin[r,s])\\
&\leq & P\left(w_n(i,\delta)\geq\eta\right)+P(\hat{\tau}_n\notin[r,s])
\end{array}
\end{equation}
where, as mentioned before $P(\hat{\tau}_n\notin[r,s])\leq\frac{\epsilon}{3}$. Furthermore, if  $\sqrt{nh_n}e_i\Delta_n(\tau h_n)$ is asymptotically stochastically uniformly equicontinuous with respect to $\tau$ on $[r,s]$, then there exists $N_{\epsilon,2}$ such that 
$$P\left(w_n(i,\delta)\geq\eta\right)\leq\frac{\epsilon}{3}$$
whenever $n\geq N_{\epsilon,2}$. Setting $N_{\epsilon}=\max\{N_{\epsilon,1},N_{\epsilon,2}\}$ we obtain that with stochastic equicontinuity we have $\sqrt{nh_n}\Delta_n(h_n)-\sqrt{nh_n}\Delta_n(\hat{\tau}_nh_n)=o_p(1)$. Now, since $\tau_1,~\tau_2,~r$ and $s$ are nonstochastic, then
$$w_n(i,\delta)=\sqrt{nh_n}\sup_{\{(\tau_1,\tau_2)\in[r,s]\times[r,s]:|\tau_1-\tau_2|<\delta \}}|e_iS_n(x;\tau_1h_n)^{-1}G_n^*(x;\tau_1h_n)-e_iS_n(x;\tau_2h_n)^{-1}G_n^*(x;\tau_2h_n)|$$
where $G_n^*(x;h_n) =(g_{n,0}^*(x;h_n),\ldots,g_{n,p}^*(x;h_n))'$. Thus, if
\begin{equation}
\label{snop1}\sup_{\{(\tau_1,\tau_2)\in[r,s]\times[r,s]:|\tau_1-\tau_2|<\delta \}}|s_{n,l}(x;\tau_1h_n)-s_{n,l}(x;\tau_2h_n)|=o_p(1)
\end{equation}
and
\begin{equation}\label{gnop1}
\sup_{\{(\tau_1,\tau_2)\in[r,s]\times[r,s]:|\tau_1-\tau_2|<\delta \}}|\sqrt{nh_n}g^*_{n,l}(x;\tau_1h_n)-\sqrt{nh_n}g^*_{n,l}(x;\tau_2h_n)|=o_p(1)
\end{equation}
the desired result is obtained.

By the Mean Value Theorem of \cite{Jennrich1969}
$$\sup_{\{(\tau_1,\tau_2)\in[r,s]^2:|\tau_1-\tau_2|<\delta\}}|s_{n,l}(x;\tau_1h_n)-s_{n,l}(x;\tau_2h_n)|\leq\sup_{\tau\in[r,s]}\left|\frac{d s_{nl}(x;\tau h_n)}{d \tau}\right|\delta~~~~~~~a.s.$$
Lemma \ref{lem:bdvar1} and Theorem 21.10 in \cite{Davidson1994} imply that equation (\ref{snop1}) holds.

Furthermore, by the Mean Value Theorem of \cite{Jennrich1969} and by Cauchy-Schwarz Inequality, we have that
\begin{equation}\begin{array}{rl}
|\sqrt{nh_n}g^*_{n,l}(x;\tau_1h_n)-\sqrt{nh_n}g^*_{n,l}(x;\tau_2h_n)|&=\left| \int_{\tau_2}^{\tau_1}\left(\sqrt{nh_n}\frac{d g^*_{nl}(x;\tau h_n)}{d\tau}\right)d\tau\right|\\
&\leq  \int_{\tau_2}^{\tau_1}\left|\sqrt{nh_n}\frac{d g^*_{nl}(x;\tau h_n)}{d\tau}\right|d\tau\\
&\leq \left|\int_{\tau_2}^{\tau_1}1d\tau\right|^{1/2}\left(  \int_{\tau_2}^{\tau_1}\left(\sqrt{nh_n}\frac{d g^*_{nl}(x;\tau h_n)}{d\tau}\right)^2d\tau\right)^{1/2}\\
&=|\tau_1-\tau_2|^{1/2}\left(  \int_{\tau_2}^{\tau_1}\left(\sqrt{nh_n}\frac{d g^*_{nl}(x;\tau h_n)}{d\tau}\right)^2d\tau\right)^{1/2}\\
&\leq |\tau_1-\tau_2|^{1/2}\left(  \int_{r}^{s}\left(\sqrt{nh_n}\frac{d g^*_{nl}(x;\tau h_n)}{d\tau}\right)^2d\tau\right)^{1/2}.
\end{array}\end{equation}
Once again, Theorem 21.10 in \cite{Davidson1994} and Lemma \ref{lem:bdvar2} imply that equation (\ref{gnop1}) holds. 
\end{proof}

\noindent\textbf{Proof of Corollary \ref{col:main}:}
\begin{proof}
For any $\epsilon>0$, given that $\hat{\tau}_n=O_p(1)$ there exists $r,s\in(0,\infty)$ with $r<s$ such that
$P(\hat{\tau}_n\notin[r,s])\leq\epsilon/3,~~~\forall~n.$.
From Theorem \ref{thm:main}, it suffices to show that
$$\sqrt{nh_n}\left(H_{n\tau}(b_n(x;\tau h_n)-b(x))-\frac{(\tau h_n)^{p+1}m^{(p+1)}(x)}{(p+1)!}+(\tau h_n)^{p+1}o_p(1) \right)$$
is stochastic equicontinuous with respect to $\tau$ on $[r,s]$ with $H_{n\tau}=diag\{(\tau h_n)^j\}_{j=0}^p$.

\cite{Masry1997}, showed that
$$H_{n\tau}(b_n(x;\tau h_n)-b(x))-\frac{(\tau h_n)^{p+1}m^{(p+1)}(x)}{(p+1)!}+(\tau h_n)^{p+1}o_p(1)=S_n^{-1}(x;\tau h_n)G_n^*(x;\tau h_n)$$
thus, from theorem \ref{thm:main}, the result follows.
\end{proof}

\clearpage

\section*{Appendix 2: Tables and Graphs}
\begin{table}[h!]
 \caption{Bias and MSE ratios for $m_1(x)$ and $m_1^{(1)}(x)$  data driven $h$ and $h_{AMISE}$}\label{tab:mc1}
\centerline{
\begin{tabular}{lclclclclclclclcl}
\\ \hline
\multirow{2}{*}{$m_1(x)$} && \multirow{2}{*}{$n$}  &&  && \multicolumn{3}{c}{$x=0.5\pi$} &&  \multicolumn{3}{c}{$x=\pi$} &&  \multicolumn{3}{c}{$x=1.5\pi$} \\ \cline{7-9} \cline{11-13} \cline{15-17}
&& && && $h_{CV}$ && $h_{ROT}$ && $h_{CV}$ && $h_{ROT}$ && $h_{CV}$ && $h_{ROT}$ \\
\hline
$\rho=0$ && 200 && 	Bias 		&& 1.014 && 0.858 && 1.267 && 1.021 && 1.022 && 0.854 \\
&&	&&				MSE 	&& 1.160 && 0.974 && 1.058 && 1.085 && 1.082 && 0.972 \\
&& 600 && 			Bias 		&& 1.015 && 0.858 && 1.196 && 1.113 && 1.015 && 0.862 \\
&&	&&				MSE 	&& 1.039 && 0.972 && 1.047 && 1.077 && 1.050 && 0.964 \\
\\
$\rho=0.5$ && 200 && Bias 		&& 1.010 && 0.856 && 1.158 && 1.115 && 1.008 && 0.859 \\
&&	&&				MSE 	&& 1.088 && 0.974 && 1.050 && 1.073 && 1.107 && 0.979 \\
&& 600 && 			Bias 		&& 0.988 && 0.861 && 1.782 && 0.807 && 1.010 && 0.856 \\
&&	&&				MSE 	&& 1.035 && 0.975 && 1.045 && 1.065 && 1.051 && 0.971 \\
\\
$\rho=0.9$ && 200 && 	Bias 		&& 0.980 && 0.794 && 1.104 && 1.033 && 0.977 && 0.835 \\
&&	&&				MSE 	&& 1.030 && 0.988 && 1.039 && 1.038 && 1.025 && 0.982 \\
&& 600 && 			Bias 		&& 0.981 && 0.851 && 0.922 && 1.013 && 0.986 && 0.872 \\
&&	&&				MSE 	&& 1.027 && 0.981 && 1.018 && 1.035 && 1.014 && 0.979 \\
\\ \hline
\multirow{2}{*}{$m_1^{(1)}(x)$} && \multirow{2}{*}{$n$}  &&  && \multicolumn{3}{c}{$x=0.5\pi$} &&  \multicolumn{3}{c}{$x=\pi$} &&  \multicolumn{3}{c}{$x=1.5\pi$} \\ \cline{7-9} \cline{11-13} \cline{15-17}
&& && && $h_{CV}$ && $h_{ROT}$ && $h_{CV}$ && $h_{ROT}$ && $h_{CV}$ && $h_{ROT}$ \\
\hline
$\rho=0$ && 200 && 	Bias 		&& 0.973 && 1.192 && 1.000 && 1.000 && 0.661 && 1.124 \\
&&	&&				MSE 	&& 1.676 && 1.237 && 1.002 && 1.002 && 1.709 && 1.253 \\
&& 600 && 			Bias 		&& 0.825 && 1.202 && 0.999 && 1.000 && 1.120 && 1.137 \\
&&	&&				MSE 	&& 1.178 && 1.240 && 1.000 && 1.001 && 1.174 && 1.237 \\
\\
$\rho=0.5$ && 200 && 	Bias 		&& 1.922 && 1.141 && 0.998 && 1.001 && 3.651 && 1.156 \\
&&	&&				MSE 	&& 1.417 && 1.250 && 0.998 && 1.003 && 1.514 && 1.245 \\
&& 600 && 			Bias 		&& 1.077 && 1.043 && 1.000 && 1.000 && 0.956 && 1.008 \\
&&	&&				MSE 	&& 1.000 && 1.001 &&1.001 && 1.001 && 1.341 && 1.242 \\
\\
$\rho=0.9$ && 200 && 	Bias 		&& 0.705 && 0.852 && 1.000 && 1.001 && 1.594 && 1.176 \\
&&	&&				MSE 	&& 1.434 && 1.272 && 1.001 && 1.004 && 1.434 && 1.290 \\
&& 600 && 			Bias 		&& 0.786 && 1.070 && 1.000 && 1.001 && 0.965 && 1.180 \\
&&	&&				MSE 	&& 1.205 && 1.248 && 1.002 && 1.002 && 1.233 && 1.256 \\
    \hline
\end{tabular}
}
\end{table}

\begin{table}[htb]
 \caption{Bias and MSE ratios for $m_2(x)$ and $m_2^{(1)}(x)$  data driven $h$ and $h_{AMISE}$}\label{tab:mc2}
\centerline{
\begin{tabular}{lclclclclclclclcl}
\\ \hline
\multirow{2}{*}{$m_2(x)$} && \multirow{2}{*}{$n$}  &&  && \multicolumn{3}{c}{$x=0.25$} &&  \multicolumn{3}{c}{$x=0.5$} &&  \multicolumn{3}{c}{$x=0.75$} \\ \cline{7-9} \cline{11-13} \cline{15-17}
&& && && $h_{CV}$ && $h_{ROT}$ && $h_{CV}$ && $h_{ROT}$ && $h_{CV}$ && $h_{ROT}$ \\
\hline
$\rho=0$ && 200 && 	Bias 		&& 1.477 && 0.902 && 0.207 && 1.194 && 1.383 && 0.924 \\
&&	&&				MSE 	&& 1.291 && 1.045 && 0.983 && 1.050 && 1.298 && 1.029 \\
 && 600 && 			Bias 		&& 1.348 && 0.962 && 0.772 && 1.414 && 1.329 && 0.957 \\
&&	&&				MSE 	&& 1.231 && 1.020 && 1.000 && 1.019 && 1.229 && 1.016 \\
\\
$\rho=0.5$ && 200 && 	Bias 		&& 1.499 && 0.916 && 0.853 && 1.058 && 1.422 && 0.922 \\
&&	&&				MSE 	&& 1.226 && 1.022 && 0.971 && 1.033 && 1.228 && 1.017 \\
 && 600 && 			Bias 		&& 1.291 && 0.966 && 0.914 && 1.025 && 1.352 && 0.963 \\
&&	&&				MSE 	&& 1.178 && 1.012 && 1.007 && 1.012 && 1.189 && 1.016 \\
\\
$\rho=0.9$ && 200 && 	Bias 		&& 1.368 && 0.877 && 0.368 && 1.163 && 1.294 && 0.898 \\
&&	&&				MSE 	&& 1.075 && 0.997 && 1.006 && 1.011 && 1.061 && 1.006 \\
 && 600 && 			Bias 		&& 1.447 && 0.941 && 0.972 && 1.017 && 1.285 && 0.955 \\
&&	&&				MSE 	&& 1.070 && 1.006 && 1.002 && 1.001 && 1.082 && 1.003 \\
    \hline
\multirow{2}{*}{$m_2^{(1)}(x)$} && \multirow{2}{*}{$n$}  &&  && \multicolumn{3}{c}{$x=0.25$} &&  \multicolumn{3}{c}{$x=0.5$} &&  \multicolumn{3}{c}{$x=0.75$} \\ \cline{7-9} \cline{11-13} \cline{15-17}
&& && && $h_{CV}$ && $h_{ROT}$ && $h_{CV}$ && $h_{ROT}$ && $h_{CV}$ && $h_{ROT}$ \\
\hline
$\rho=0$ && 200 && 	Bias 		&& 2.119 && 0.943 && 1.607 && 0.878 && 2.909 && 0.796 \\
&&	&&				MSE 	&& 1.514 && 1.133 && 1.871 && 1.133 && 1.481 && 1.163 \\
 && 600 && 			Bias 		&& 1.619 && 1.018 && 1.562 && 0.969 && 1.269 && 0.937 \\
&&	&&				MSE 	&& 1.751 && 1.067 && 2.413 && 1.067 && 1.422 && 1.050 \\
\\
$\rho=0.5$ && 200 && 	Bias 		&& 1.808 && 0.836 && 1.517 && 0.950 && 3.076 && 0.908 \\
&&	&&				MSE 	&& 1.765 && 1.120 && 1.643 && 1.120 && 1.573 && 1.122 \\
 && 600 && 			Bias 		&& 2.272 && 0.918 && 1.172 && 0.940 && 3.139 && 0.832 \\
&&	&&				MSE 	&& 1.492 && 1.052 && 1.715 && 1.052 && 1.667 && 1.050 \\
\\
$\rho=0.9$ && 200 && 	Bias 		&& 1.681 && 0.915 && 1.538 && 0.904 && 2.035 && 0.870 \\
&&	&&				MSE 	&& 1.972 && 1.207 && 2.695 && 1.207 && 1.629 && 1.182 \\
 && 600 && 			Bias 		&& 1.622 && 0.948 && 1.206 && 0.969 && 2.285 && 1.007 \\
&&	&&				MSE 	&& 2.103 && 1.066 && 1.486 && 1.066 && 1.691 && 1.076 \\
    \hline
\end{tabular}
}
\end{table}

\clearpage

\clearpage
\setlength{\baselineskip}{12pt}
\bibliographystyle{elsart-harv.bst}
\bibliography{library}

\begin{thebibliography}{18}
\expandafter\ifx\csname natexlab\endcsname\relax\def\natexlab#1{#1}\fi
\expandafter\ifx\csname url\endcsname\relax
  \def\url#1{\texttt{#1}}\fi
\expandafter\ifx\csname urlprefix\endcsname\relax\def\urlprefix{URL }\fi

\bibitem[{Bernstein(1927)}]{Bernstein1927}
Bernstein, S., 1927. {Sur l'extension du theorem du calcul des probabilites aux
  sommes de quantites dependantes}. Mathematische Annalen 97, 1--59.

\bibitem[{Boente and Fraiman(1995)}]{Boente1995}
Boente, G., Fraiman, R., 1995. {Asymptotic distribution of data-driven
  smoothers in density and regression estimation under dependence}. Canadian
  Journal of Statistics 23, 383--397.

\bibitem[{Bosq(1998)}]{Bosq1998}
Bosq, D., 1998. {Nonparametric statistics for stochastic processes: estimation
  and prediction}. No. 110 in Lecture Notes in Statistics. Springer Verlag.

\bibitem[{Davidson(1994)}]{Davidson1994}
Davidson, J., 1994. {Stochastic Limit Theory}. Oxford University Press, New
  York.

\bibitem[{Dony et~al.(2006)Dony, Einmahl, and Mason}]{Dony2006}
Dony, J., Einmahl, U., Mason, D., 2006. {Uniform in bandwidth consistency of
  local polynomial regression function estimators}. Austrian Journal of
  Statistics 35, 105--120.

\bibitem[{Doukhan(1994)}]{Doukhan1994}
Doukhan, P., 1994. {Mixing}. Springer-Verlag, New York.

\bibitem[{Fan(1992)}]{Fan1992}
Fan, J., 1992. {Design-adaptive nonparametric regression}. Journal of the
  American Satistical Association 87, 998--1004.

\bibitem[{Jennrich(1969)}]{Jennrich1969}
Jennrich, R.~I., 1969. {Asymptotic Properties of Non-Linear Least Squares
  Estimators}. The Annals of Mathematical Statistics 40~(2), 633--643.

\bibitem[{Li and Racine(2004)}]{Li2004}
Li, Q., Racine, J., 2004. {Cross-validated local linear nonparametric
  regression}. Statistica Sinica 14, 485--512.

\bibitem[{Martins-Filho and Yao(2009)}]{Martinsfilho2009}
Martins-Filho, C., Yao, F., 2009. {Nonparametric regression estimation with
  general parametric error covariance}. Journal of Multivariate Analysis 100,
  309--333.

\bibitem[{Masry and Fan(1997)}]{Masry1997}
Masry, E., Fan, J., 1997. {Local polynomial estimation of regression functions
  for mixing processes}. Scandinavian Journal of Statistics 24, 1965--1979.

\bibitem[{Parzen(1962)}]{Parzen1962}
Parzen, E., 1962. {On estimation of a probability density and mode}. Annals of
  Mathematical Statistics 33, 1065--1076.

\bibitem[{Pham and Tran(1985)}]{Pham1985}
Pham, T.~D., Tran, L.~T., 1985. {Some mixing properties of time series models}.
  Stochastic Processes and their Applications 19, 297--303.

\bibitem[{Robinson(1983)}]{Robinson1983}
Robinson, P.~M., 1983. {Nonparametric estimators for time series}. Journal of
  Time Series Analysis 4, 185--207.

\bibitem[{Ruppert et~al.(1995)Ruppert, Sheather, and Wand}]{Ruppert1995}
Ruppert, D., Sheather, S., Wand, M.~P., 1995. {An effective bandwidth selector
  for local least squares regression}. Journal of the American Statistical
  Association 90, 1257--1270.

\bibitem[{Silverman(1986)}]{Silverman1986}
Silverman, B.~W., 1986. {Density estimation for statistics and data analysis}.
  Chapman and Hall, London.

\bibitem[{Xia and Li(2002)}]{Xia2002}
Xia, Y., Li, W.~K., 2002. {Asymptotic Behavior of Bandwidth Selected by the
  Cross-Validation Method for Local Polynomial Fitting}. Journal of
  Multivariate Analysis 83, 265--287.

\bibitem[{Ziegler(2004)}]{Ziegler2004}
Ziegler, K., 2004. {Adaptive kernel estimation of the mode in nonparametric
  random design regression model}. Probability and Mathematical Statistics 24,
  213--235.

\end{thebibliography}
\end{document}